\documentclass[12pt]{article}
\usepackage[T2A]{fontenc}
\usepackage[utf8]{inputenc}
\usepackage{amssymb,amsfonts,amsmath,color}
\usepackage{amssymb,amsmath,amsfonts,amsthm,amscd,latexsym,indentfirst,verbatim,mathtext,cite,enumerate,float}
\def\RB{\textrm{RB}}
\def\rb{\textrm{rb}}
\def\id{\textrm{id}}
\def\diag{\textrm{diag}}
\def\charr{\textrm{char}\,}
\def\Var{\mathrm{Var}}
\def\Span{\mathrm{Span}}
\def\Aut{\mathrm{Aut}}
\def\Imm{\mathrm{Im}\,}

\def\Spec{\mathrm{Spec}\,}

\begin{document}

\begin{center}
{\Large
Spectrum of Rota---Baxter operators}

V. Gubarev
\end{center}

\begin{abstract}
Rota---Baxter operators have been studied since 1960, and there are lots of
their applications in mathematical physics, number theory, and noncommutative geometry.
We state a surprisingly general property of such operators:
the spectrum of every Rota---Baxter operator of weight~$\lambda$ on a~finite-dimensional unital
(not necessarily associative) algebra is a~subset of $\{0,-\lambda\}$.
We even extend this result on the case of infinite-dimensional algebraic algebras (when $\charr F = 0$).
Based on these results, we define a new invariant of an algebra: the Rota---Baxter $\lambda$-index
$\rb_\lambda(A)$ of an algebra $A$ as the infimum of the degrees
of minimal polynomials of all Rota---Baxter operators of weight~$\lambda$ on~$A$.
We calculate the Rota---Baxter $\lambda$-index for the matrix algebra $M_n(F)$, $\charr F = 0$:
it is shown that $\rb_\lambda(M_n(F)) = 2n-1$.

\end{abstract}

\section{Introduction}

Given an algebra $A$ over a ground field~$F$ and a scalar $\lambda\in F$,
a~linear operator $R\colon A\rightarrow A$ is called a Rota---Baxter operator
(RB-operator, for short) of weight~$\lambda$ on $A$ if the following identity holds for all $x,y\in A$
\begin{equation}\label{RB}
R(x)R(y) = R( R(x)y + xR(y) + \lambda xy ).
\end{equation}
An algebra $A$ endowed with an RB-operator is called Rota---Baxter algebra (RB-algebra).

Rota---Baxter operators appeared in the work of G.~Baxter~\cite{Baxter} in 1960
in the study of fluctuation theory as a formal generalization of integration by parts formula.
Further, they were studied by G.-C.~Rota~\cite{Rota} and others \cite{Atkinson,Miller66,Cartier}.
The deep connections between solutions of different versions
of the classical Yang---Baxter equation in mathematical physics
and RB-operators were discovered in~\cite{BelaDrin82,Aguiar00,FardThesis}.
The  subject of double algebras is also tightly connected to Rota---Baxter operators,
see \cite{DoublePoissonFree,DoubleLie}.
Since 2000, the connection between RB-algebras and so called Loday algebras
(prealgebras and postalgebras) has been studied~\cite{Aguiar00,GolubchikSokolov,BBGN2011,Embedding}.

Rota---Baxter operators have a broad area of applications
in symmetric polynomials, shuffle algebra, etc.~\cite{Atkinson,GuoMonograph,GuoKeigher}.

There is a plenty of papers devoted to a partial or complete
classification of Rota---Baxter operators on a given small-dimensional algebra:
$\mathrm{sl}_2(F)$~\cite{Kolesnikov,sl2,sl2-0},
$\mathrm{sl}_3(F)$~\cite{KonovDissert},
$M_2(F)$~\cite{Aguiar00-2,BGP,Mat2},
$M_3(F)$~\cite{GonGub,Sokolov},
the Grassmann algebra $\mathrm{Gr}_2$~\cite{BGP},
and others, see the survey~\cite{Unital} and the references therein.

The classifications obtained in these works are often bulky and hard to observe,
so there is a need to investigate the general properties of Rota---Baxter operators
at least on finite-dimensional (and close to them) algebras.
For this, the author proved in~\cite{Unital} the following result.
Given a unital algebraic associative (alternative, Jordan) algebra $A$
over a~field of characteristic zero,
every RB-operator of weight zero on $A$ is nilpotent.
Moreover, if $A$ is algebraic of a bounded degree, then
there exists $N$ such that $R^N = 0$ for every RB-operator~$R$ of weight zero on $A$.
(By an algebraic algebra in the nonassociative case we mean such an algebra
that its every finite subset generates a~finite-dimensional subalgebra.)

The statement above leads to a new invariant of an algebra~$A$:
the Rota---Baxter index (RB-index) $\rb(A)$ defined
as the upper bound of nilpotency indices of all RB-operators of weight zero on~$A$.
This invariant was calculated for different finite-dimensional
algebras~\cite{Unital,Gub2017}, including the matrix algebra $M_n(F)$.

Now, we prove that the spectrum of a Rota---Baxter of nonzero weight~$\lambda$ on
a finite-dimensional unital (or algebraic unital over a field of characteristic zero) algebra
is a~subset of $\{0,-\lambda\}$. The technique of the proof is quite different from the previous one.
Moreover, we avoid the two conditions from the result from~\cite{Unital} for the case~$\lambda = 0$.
Firstly, we remove the restriction on a variety of algebras, i.\,e., we do not assume
that $A$~is associative, alternative, or Jordan.
Secondly, we state the result in the case of a~field of any characteristic.
Note that the proofs in the case of positive characteristic is very unlike the proof when $\charr F = 0$.

Before stating the main result of the work, we need a new definition.
Given an algebra $A$, denote the set of all RB-operators of weight~$\lambda$ on~$A$ by $\RB_\lambda(A)$.
Define the Rota---Baxter $\lambda$-index ($\RB(\lambda)$-index) of $A$ as
$$
\rb_\lambda(A) = \min\{n\in\mathbb{N}\mid \mbox{for all }R\in \RB_\lambda(A)\
 \mbox{exists } k:\ R^k(R+\lambda\id)^{n-k} = 0\}.
$$
If such a number does not exist, then put $\rb_\lambda(A) = \infty$.

If $\lambda = 0$, then $\rb_0(A)$ coincides with
the earlier defined Rota---Baxter index $\rb(A)$~\cite{Unital}.

In~\S3 (Theorem~1, over a~field of characteristic zero)
and in~\S4 (Theorem~2, over a~field of positive characteristic),
we prove

{\bf Main Theorem}.
Let $A$ be a finite-dimensional unital algebra over a field~$F$. Let us fix $\lambda\in F$.

a) Given an RB-operator $R$ of weight~$\lambda$ on $A$,
there exist $k,l\geq0$ such that $R^k(R+\lambda\id)^l = 0$. In particular,
$\Spec(R)\subset\{0,-\lambda\}$.

b) There exists a natural~$n$ such that for every RB-operator~$R$
of weight~$\lambda$ on $A$ one can find $0\leq k\leq n$ such that
$R^k(R+\lambda\id)^{n-k}(1) = 0$ holds.
Moreover, $\rb_\lambda(A)\leq 2n$.

When $\charr F = 0$, the part a) of the Theorem holds even if $A$ is algebraic (see~\S3).

In~\S5 (Theorem~3), we state another important property of Rota---Baxter operators.
Given a unital algebraic algebra $A$ and an RB-operator~$R$ of nonzero weight on~$A$,
\underline{the minimal polynomial of the element $R(1)$ has no multiple roots}.
This result allows us to rethink two of the three main results from~\cite{Unital}.
One of them said that every Rota---Baxter operator of nonzero weight on the Grassmann algebra is splitting.
Another one was diagonalizibilty of $R(1)$ for an RB-operator~$R$ of nonzero weight
on the matrix algebra $M_n(F)$ over a field of characteristic zero.
Now, we clarify the common nature of both these results, which concerns
the minimal polynomial of $R(1)$.

In~\S6, we compute the Rota---Baxter $\lambda$-index for different algebras:
a sum of fields, the matrix algebra $M_n(F)$, its subalgebra $t_n(F)$ of not strictly upper-triangular matrices,
the simple Jordan algebra $J_n(f)$ of a nondegenerate symmetric bilinear form~$f$,
and the matrix (split) Cayley---Dickson algebra $C(F)$.

Let us collect all results obtained earlier (see~\cite{BGP,Gub2017,Unital}) and
in the current work about $\RB(\lambda)$-index of particular algebras in the following table (see the next page).
By $\mathrm{Gr}_n$ we denote the Grassmann algebra
of a vector space $V = \Span\{e_1,\ldots,e_n\}$.
The sign ``(?)'' denotes the hypothesis,
see Conjectures~1 and~2 in~\S6.2.

In~\S7, we apply the Main Theorem for Rota---Baxter modules (\S7.1)
and parameterized Rota---Baxter operators (\S7.2).

In~\S7.3, we generalize the construction which is already known
for different varieties of algebras~\cite{GuoMonograph,preLie,Burde}.
Given an algebra~$A$ of a~variety $\Var$ and an RB-operator~$R$ on~$A$,
one get a countable sequence of algebras~$A_i$ belonging to the same variety~$\Var$.
Moreover, all of $A_i$ considered with~$R$ are Rota---Baxter algebras.

We also apply the Main Theorem to state the structure of the limiting algebra $A_\infty$
(see~\S7.3).

\begin{table}[t]
\begin{center}
Table: the Rota---Baxter $\lambda$-index of different algebras
\end{center}
\begin{tabular}{c|c|c|c}
Algebra & Field & $\rb_0$ & $\rb_\lambda$, $\lambda\neq0$ \\
\hline
$F^n$ & any & 1 & $n$ \\
$M_n(F)$ & $\charr(F) = 0$ & $2n-1$ & $2n-1$ \\
$t_n(F)$ & $\charr(F) = 0$ & $n$ (?) & $2n-1$ \\
$\mathrm{Gr}_n$ & $\charr(F) = 0$ & ? & 2 \\
$J_n(f)$ & $\bar{F} {=} F$, $\charr(F) {\neq}2,3$ & $\begin{cases}
2, & n=3  \\
3, & n>3 \end{cases}$ & $\begin{cases}
2, & 2\!\not|\,n  \\
3, & 2\,|\,n
\end{cases}$\\
$C(F)$ & $\bar{F} {=} F$, $\charr(F) {\neq}2,3$ & 3 & 3 \\
$A = \bigoplus\limits_{i=1}^k M_{n_i}(F)$ & $\charr(F) = 0$ &
$2\max\limits_{i=1,\ldots,k}\{n_i\} {-} 1$\,(?) &
$\sum\limits_{i=1}^k n_i {+} \max\limits_{i=1,\ldots,k}\{n_i\} - 1$\,(?)
\end{tabular}
\end{table}

At the end of Introduction, let us clarify other applications of the stated spectral results about Rota---Baxter operators on unital algebras including the results obtained after the work appeared on arxiv.org.

The knowledge of the spectrum of an RB-operator on a given finite-dimensional unital algebra~$A$ really helps to classify all RB-operators on it.
As a consequence, we may list all possible Jordan forms of RB-operators and then consider such forms case-by-case. 

The spectral results have been already applied for the classification of RB-operators on the following algebras,
\begin{itemize}
 \item the simple Jordan superalgebra $D_t$~\cite{Dt},
 \item Jordan algebras of dimension~2 and~3~\cite{Sun}.
\end{itemize}

Another conclusion of the obtained theory is the following: RB-operators on which the maximal value of the Rota---Baxter $\lambda$-index on~$A$ is reached are of special interest. 
In~\cite{Gub2021,GonGub2}, Rota---Baxter operators $R_n(\lambda)$ of weight~$\lambda$ on the matrix algebra~$M_n(F)$ providing the maximal value
$\rb_\lambda(M_n(F)) = 2n-1$ were put into the machinery connecting them and double Lie algebras. 
As an outcome, for $\lambda = 0$, the first simple double Lie algebra was found. 
For $\lambda\neq0$, the example of almost simple $\lambda$-double Lie algebra was constructed.
The Main Theorem was applied in~\cite{GonGub2} to prove that there are no simple finite-dimensional $\lambda$-double Lie algebras for nonzero~$\lambda$.

\section{Preliminaries}

Trivial RB-operators of weight $\lambda$ are zero operator and $-\lambda\id$.

{\bf Lemma 1}~\cite{GuoMonograph}.
Given an RB-operator $P$ of weight $\lambda$,

(a) the operator $-P-\lambda\id$ is an RB-operator of weight $\lambda$,

(b) the operator $\lambda^{-1}P$ is an RB-operator of weight 1, provided $\lambda\neq0$.

Given an algebra $A$, let us define a map $\phi$ on the set of all RB-operators on $A$
as $\phi(P)=-P-\lambda(P)\id$, where $\lambda(P)$ denotes the weight of an RB-operator~$P$.
It is clear that $\phi^2$ coincides with the identity map.

{\bf Lemma 2}~\cite{BGP}.
Given an algebra $A$, an RB-operator $P$ on $A$ of weight $\lambda$,
and $\psi\in\Aut(A)$, the operator $P^{(\psi)} = \psi^{-1}P\psi$
is an RB-operator of weight~$\lambda$ on~$A$.

The same result is true when $\psi$ is an antiautomorphism of $A$, i.e.,
a~bijection from $A$ to $A$ satisfying $\psi(xy) = \psi(y)\psi(x)$ for all $x,y\in A$;
e.g., transpose on the matrix algebra.

{\bf Lemma 3}~\cite{GuoMonograph}.
Let an algebra $A$ split as a vector space
into a direct sum of two subalgebras $A_1$ and $A_2$.
An operator $P$ defined as
\begin{equation}\label{Split}
P(a_1 + a_2) = -\lambda a_2,\quad a_1\in A_1,\ a_2\in A_2,
\end{equation}
is RB-operator of weight~$\lambda$ on~$A$.

Let us call an RB-operator from Lemma~3 as
{\it splitting} RB-operator with subalgebras $A_1,A_2$.
Note that the set of all splitting RB-operators on
an algebra $A$ is in bijection with all decompositions of~$A$
into a~direct sum of two subalgebras $A_1,A_2$.

{\bf Lemma 4}~\cite{BGP}.
Let $A$ be a unital algebra and let $R$ be an RB-operator
of nonzero weight~$\lambda$ on $A$.
If $R(1)\in F$, then $R$ is splitting with one of subalgebras containing~1.

The following construction of RB-operators of nonzero weight generalizes Lemma~3.

{\bf Lemma 5}~\cite{Guil}.
Let an algebra $A$ be a direct sum of subspaces $A_-,A_0,A_+$,
let $A_\pm,A_0$ be subalgebras of $A$,
and let $A_\pm$ be $A_0$-modules.
If $R_0$ is an RB-operator of weight~$\lambda$ on $A_0$,
then an operator $P$ defined as
\begin{equation}\label{RB:SubAlg2}
P(a_-+a_0+a_+) = R_0(a_0) - \lambda a_+,\quad
a_{\pm}\in A_{\pm},\ a_0\in A_0,
\end{equation}
is an RB-operator of weight~$\lambda$ on~$A$.

Let us call an RB-operator of nonzero weight
defined by~\eqref{RB:SubAlg2} as triangular-splitting one
provided that at least one of $A_-,A_+$ is nonzero.

{\bf Lemma 6}.
Let $A$ be a unital algebra and $R$ be an RB-operator of weight~$\lambda$ on~$A$.

a) Then the subalgebra $S$ generated by the set $\{R^k(1) \mid k\in\mathbb{N}\}$
is commutative and associative.

b) \cite[\S3.2.2]{GuoMonograph} for all $k,l\in\mathbb{N}$
the following formula holds,
\begin{equation}\label{R^k(1)R^l(1)Main}
R^k(1)R^l(1)
 =  \sum\limits_{j=0}^{\min\{k,l\}}\lambda^j\binom{k}{j}\binom{k+l-j}{k}R^{k+l-j}(1).
\end{equation}

{\sc Proof}.
a) We prove that the associator
$$
h(k,l,m)
 = (R^k(1),R^l(1),R^m(1))
 = (R^k(1)R^l(1))R^m(1) - R^k(1)(R^l(1)R^m(1))
$$
is zero for all natural $k,l,m$ by induction on $r = k+l+m$.
For $r = 0$, it is true. Suppose that the statement holds true for all $p<r$.
Note that $h(k,l,m) = 0$ if at least one of the numbers $k,l,m$ is zero.
So, we may assume that $k,l,m\geq1$.

On the one hand, we have
\begin{multline*}
(R^k(1)R^l(1))R^m(1) \\
 = R( R^{k-1}(1)R^l(1) + R^k(1)R^{l-1}(1) + \lambda R^{k-1}(1)R^{l-1}(1))R^m(1) \\
 = R[ (R^{k-1}(1)R^l(1))R^m(1) + (R^k(1)R^{l-1}(1))R^m(1) + \lambda (R^{k-1}(1)R^{l-1}(1))R^m(1) \\
  + R( R^{k-1}(1)R^l(1) + R^k(1)R^{l-1}(1) + \lambda R^{k-1}(1)R^{l-1}(1))R^{m-1}(1) \\
  + \lambda (R^{k-1}(1)R^l(1))R^{m-1}(1)
  + \lambda(R^k(1)R^{l-1}(1))R^{m-1}(1)
  + \lambda^2 (R^{k-1}(1)R^{l-1}(1))R^{m-1}(1) ].
\end{multline*}
On the other hand,
\begin{multline*}
R^k(1)(R^l(1)R^m(1)) \\
 = R^k(1) R( R^{l-1}(1)R^m(1) + R^l(1)R^{m-1}(1) + \lambda R^{l-1}(1)R^{m-1}(1)) \\
 = R[ R^k(1)(R^{l-1}(1)R^m(1)) + R^k(1)(R^l(1)R^{m-1}(1)) + \lambda R^k(1)(R^{l-1}(1)R^{m-1}(1)) \\
  + R^{k-1}(1) R(  R^{l-1}(1)R^m(1) + R^l(1)R^{m-1}(1) + \lambda R^{l-1}(1)R^{m-1}(1) ) \\
  + \lambda R^{k-1}(1)(R^{l-1}(1)R^m(1))
  + \lambda R^{k-1}(1)(R^l(1)R^{m-1}(1))
  + \lambda^2 R^{k-1}(1)(R^{l-1}(1)R^{m-1}(1)) ].
\end{multline*}
Applying the induction hypothesis and~\eqref{RB},
we compute $h(k,l,m) = R(\Delta)$, where
\begin{multline*}
\Delta
 = R^{k-1}(1)R^l(1)R^m(1) + R^k(1)R^{l-1}(1)R^m(1) + \lambda R^{k-1}(1)R^{l-1}(1)R^m(1) \\
 + R^{k}(1)R^l(1)R^{m-1}(1) + \lambda R^{k-1}(1)R^l(1)R^{m-1}(1) \\
 + \lambda R^k(1)R^{l-1}(1)R^{m-1}(1) + \lambda^2 R^{k-1}(1)R^{l-1}(1)R^{m-1}(1) \\ \allowdisplaybreaks
 - R^k(1)R^{l-1}(1)R^m(1) - R^k(1)R^l(1)R^{m-1}(1) - \lambda R^k(1)R^{l-1}(1)R^{m-1}(1) \\
 - R^{k-1}(1)R^l(1)R^{m}(1) - \lambda R^{k-1}(1)R^{l-1}(1)R^m(1) \\
 - \lambda R^{k-1}(1)R^l(1)R^{m-1}(1)
  - \lambda^2 R^{k-1}(1)R^{l-1}(1)R^{m-1}(1)
 = 0.
\end{multline*}

Finally, we prove that $R^k(1)R^l(1) = R^l(1)R^k(1)$ by induction on $r = k+l$.
For $r=0$, it is trivial. Suppose commutativity holds for all $p<r$.
If $k = 0$ or $l = 0$, then $R^k(1)R^l(1) = R^l(1)R^k(1)$.
Assume that $k,l\geq1$. We have
\begin{multline*}
R^k(1)R^l(1) - R^l(1)R^k(1) \\
 = R( R^{k-1}(1)R^l(1) + R^k(1)R^{l-1}(1) + \lambda R^{k-1}(1)R^{l-1}(1) \\
 - R^{l-1}(1)R^k(1) - R^l(1)R^{k-1}(1) - \lambda R^{l-1}(1)R^{k-1}(1) )
 = 0
\end{multline*}
by the induction hypothesis. \hfill $\square$

\section{Main Theorem in characteristic zero}

Let $A$ be a unital associative algebra and $R$ be an RB-operator of weight~$\lambda$ on~$A$.
It is known~\cite{GuoMonograph} that the following formula
\begin{equation}\label{Stirling1}
n! R^n(1)  = \sum\limits_{k=1}^n (-1)^{n-k}\lambda^{n-k}s(n,k)(R(1))^k
\end{equation}
holds. Here $s(n,k)$ is a~Stirling number of the first kind.
By Lemma~6, we may remove the restriction on~$A$ to be associative.

{\bf Lemma 7}.
Let $A$ be a unital algebraic algebra over a field of characteristic zero
and let $R$ be an RB-operator of a nonzero weight~$\lambda$ on $A$.
Then there exist $r,s$ such that $(R+\lambda\id)^r R^s(1) = 0$.

{\sc Proof}.
Denote $a = R(1)$. We put $\lambda = -1$ for simplicity.
Since $A$ is algebraic, there exists
a minimal polynomial $m$ for $a$ of the degree $n$.
The formula~\eqref{Stirling1} allows us
to present~$m$ in the form
\begin{equation}\label{MinPolWeightNonzero}
m(a) = R^n(1) + \alpha_{n-1}R^{n-1}(1) + \ldots + \alpha_1 R^1(1) + \alpha_0 = 0,\quad \alpha_i\in F.
\end{equation}
Applying the formula
\begin{equation}\label{R^k(1)R(1)}
R^k(1)R(1) = (k+1)R^{k+1}(1) + \lambda k R^k(1)
\end{equation}
holding by Lemma~6b, we calculate
\begin{multline}\label{MinPolWeightNonzeroXR(1)}
0 = m(a)R(1)
 = (n+1)R^{n+1}(1) + n(\alpha_{n-1}-1)R^n(1) + (n-1)(\alpha_{n-2}-\alpha_{n-1})R^{n-1}(1) \\
 +\ldots + \alpha_0 R^1(1).
\end{multline}
Subtracting~\eqref{MinPolWeightNonzeroXR(1)} from $(n+1)R(m(a))$, we obtain
\begin{equation}\label{MinPolWeightNonzeroXR(1)Reduced}
0 = (\alpha_{n-1}+n)R^n(1) + (2\alpha_{n-2}+(n-1)\alpha_{n-1})R^{n-1}(1) + \ldots + (n\alpha_0 + \alpha_1)R^1(1).
\end{equation}
We have two cases.

{\sc Case 1}: all coefficients of~\eqref{MinPolWeightNonzeroXR(1)Reduced} are zero.
Then it is easy to derive the formulas
$$
\alpha_{n-1} = - n,\quad
\alpha_{n-2} = \frac{n(n-1)}{2},\quad \ldots,\quad
\alpha_{n-i} = (-1)^i \binom{n}{i},\quad \ldots,\quad
\alpha_0 = (-1)^n.
$$
So, we get the equality $(R-\id)^n(1) = 0$, and we are done.

{\sc Case 2}: not all coefficients of~\eqref{MinPolWeightNonzeroXR(1)Reduced} are zero.
By minimality of $m$, it means that $m$ and~\eqref{MinPolWeightNonzeroXR(1)Reduced} are proportional.
Find a minimal $k$ such that $\alpha_0 = \ldots = \alpha_{k-1} = 0$ and $\alpha_k\neq0$.
Thus, coefficients by $R^i(1)$ in $km(a)$ and~\eqref{MinPolWeightNonzeroXR(1)Reduced} are equal.
Again, we write down the equalities
\begin{gather*}
\alpha_{k+1} = -(n-k)\alpha_k,\quad
\alpha_{k+2} = \frac{(n-k)(n-k-1)}{2}\alpha_k,\quad \ldots,\\
\alpha_{n-i} = (-1)^{n-k-i}\binom{n-k}{i}\alpha_k,\quad \ldots,\quad
\alpha_{n-1} = (-1)^{n-k-1}(n-k)\alpha_k.
\end{gather*}
Thus, we get
$\frac{(-1)^{n-k}}{\alpha_k}m(a) = (R-\id)^{n-k}R^k(1) = 0$. \hfill $\square$

\newpage

{\bf Lemma 8}.
Let $A$ be an algebra and $R$ be an RB-operator of weight~$\lambda$ on $A$.

a) If $\lambda = 0$, then $\ker(R^k)$ is an $\Imm(R^k)$-module for each $k\geq1$.

b)~\cite{Burde} If $\lambda \neq 0$, then $\ker(R^k)$ is an ideal in $\Imm(R+\lambda\id)^k$.

{\sc Proof}.
We prove Lemma~8 by induction on $k$. The following equality
\begin{equation}\label{kerRnIdeal}
R(x(R(y)+\lambda y))  = R(x)R(y) - R(R(x)y)
\end{equation}
holds by~\eqref{RB} for all $x,y\in A$.
If $x\in \ker(R)$, then $x\Imm(R+\lambda\id)\in \ker(R)$.
Analogously, $\Imm(R+\lambda\id)x\in \ker(R)$.

Suppose we have proved Lemma~8 for all numbers less or equal to~$k$.
Let $x\in\ker(R^{k+1})$, $y\in A$ and $z = R(y)$.
Then by~\eqref{kerRnIdeal} we get
$$
R^{k+1}(x(R+\lambda\id)^{k+1}(y))
 = R^k[ R(x)(R+\lambda\id)^k(z) - R(R(x)(R+\lambda\id)^k(y)) ].
$$
By the induction hypothesis,
$x\Imm(R+\lambda\id)^{k+1}\in \ker(R^{k+1})$.
Analogously, we obtain $\Imm(R+\lambda\id)^{k+1}x\in \ker(R^{k+1})$. \hfill $\square$

Given an algebra $A$, denote the set of all RB-operators of weight~$\lambda$ on~$A$ by $\RB_\lambda(A)$.
Define Rota---Baxter $\lambda$-index (RB($\lambda$)-index) of $A$ as
$$
\rb_\lambda(A) = \min\{n\in\mathbb{N}\mid \mbox{for all }R\in \RB_\lambda(A)\
 \mbox{exists } k:\ R^k(R+\lambda\id)^{n-k} = 0\}.
$$
If such a number does not exist, then put $\rb_\lambda(A) = \infty$.

When $\lambda = 0$,
$$
\rb_0(A) = \min\{n\in\mathbb{N}\mid  R^n = 0\mbox{ for all }R\in \RB_0(A)\}
$$
which coincides with the earlier defined Rota---Baxter index of $A$~\cite{Unital}
(it was denoted there as $\rb(A)$).

Now, we prove the Main Theorem over a field of characteristic zero, generalizing~\cite{Unital}.

{\bf Theorem 1}.
Let $A$ be a unital algebraic algebra over
a field $F$ of characteristic zero. Let us fix $\lambda\in F$.

a) Given an RB-operator $R$ of weight~$\lambda$ on $A$,
there exist $k,l\geq0$ such that $R^k(R+\lambda\id)^l = 0$. In particular,
$\Spec(R)\subset\{0,-\lambda\}$.

b) Suppose that $A$ is algebraic of the restricted degree $N$, then
there exists a natural $n\leq N$ such that for every RB-operator
$R$ of weight~$\lambda$ on $A$ one can find $0\leq k\leq n$
such that
$R^k(R+\lambda\id)^{n-k}(1) = 0$ holds.
Moreover, $\rb_\lambda(A)\leq 2n$.

{\sc Proof}.
Let $\lambda = 0$. In~\cite{Unital}, it was proved that
$R^n(1) = 0$  for some $n\in\mathbb{N}$ under the condition that $A$ is power-associative.
By Lemma~6, we may avoid this condition.
When $A$ is algebraic of the restricted degree~$N$, then $n\leq N$.
Let $x\in \Imm(R^n)$. By Lemma~8a, $1\cdot x = x \in \ker(R^n)$.
Thus, $R^{2n}(A) = (0)$, we are done.

Let $\lambda\neq0$. We may assume that $\lambda = -1$.
Since $A$ is algebraic, we may consider the minimal (unital)
polynomial $f(x)$ of the element $R(1) = a$ of the degree $n$.
By Lemma~7, there exists $k$ such that
$R^k(R-\id)^{n-k}(1) = 0$.
When $A$ is algebraic of the restricted degree~$N$, then $n\leq N$.

We want to prove that $1\in \ker(R^k)\oplus \ker(R-\id)^{n-k}$.
First, it is easy to check that $\ker(R^k)\cap \ker(R-\id)^{n-k} = (0)$.
Second, let $e_\gamma$, $\gamma\in\Gamma$, and $f_\delta$, $\delta\in\Delta$,
be a~linear basis of $\ker(R^k)$ and of $\ker(R-\id)^{n-k}$ respectively.
We may complete the set
$\{e_\gamma\mid \gamma\in\Gamma\}\cup \{f_\delta\mid \delta\in\Delta\}$
to a basis of $A$ by adding vectors $\{g_\pi\mid \pi\in\Pi\}$.
So, we represent 1 as
$$
1 = 1_e + 1_f + 1_g
$$
meaning that $1_e\in \Span\{e_\gamma\}$, $1_f\in \Span\{f_\delta\}$,
and $1_g\in \Span\{g_\pi\}$. Let us show that $1_g = 0$ and so,
$1\in \ker(R^k)\oplus \ker(R-\id)^{n-k}$. Indeed,
$R^k(1) = R^k(1_f) + R^k(1_g)$, and $R^k(1_f)\in \ker(R-\id)^{n-k}$,
since $\ker(R-\id)^{n-k}$ is $R$-invariant.
As $(R-\id)^{n-k}(R^k(1)) = 0$, we conclude that
$R^k(1_g) \in \ker(R-\id)^{n-k}$.
The map $R^k$ is invertible on $\ker(R-\id)^{n-k}$, so, $1_g\in \ker(R-\id)^{n-k}$
and $1_g = 0$.

Finally, let $x\in \Imm R^{n-k}(R-\id)^k$.
By Lemma~8b,
$$
(\ker(R^k)\oplus \ker(R-\id)^{n-k})\cdot\Imm R^{n-k}(R-\id)^k
 \subset \ker(R^k)\oplus \ker(R-\id)^{n-k}.
$$
So, $1\cdot x = x \in \ker(R^k)\oplus \ker(R-\id)^{n-k}$.
As a consequence, $R^n(R-\id)^n(A) = (0)$.
\hfill $\square$

{\bf Remark 1}.
Let us show that all conditions of Theorem~1 are necessary.

\begin{itemize}
\item
If $A$ is not unital, then we may consider an algebra $A$ with zero product.
Then every linear map on $A$ is an RB-operator of a weight~$\lambda$.
So, we may find a linear map satisfying $\Spec(R)\not\subset\{0,-\lambda\}$.

\item
If $A$ is not algebraic, then $A = F[x]$ with the RB-operator $R(x^n) = \frac{x^{n+1}}{n+1}$
is a~counterexample to the conclusion of Theorem~1 when $\lambda = 0$.
Regardless of weight~$\lambda$, we may consider $A$ equal
to the free commutative unital RB-algebra generated by
a non-empty set $X = \{1\}$ (see~\cite{GuoKeigher}),
it is also a~counterexample to the conclusion of Theorem~1.

\item
If $\charr (F) = p > 0$, then let $A$ be
the free commutative unital RB-algebra of weight~0 generated by unit.
Due to~\cite{GuoKeigher}, $A = \bigoplus\limits_{i\geq0} Fe_i$ with the product
$e_i e_j = \left(\binom{i+j}{i}\!\! \mod p\right)e_{i+j}$
and the RB-operator~$R$ acting as follows, $R(e_i) = e_{i+1}$.
Note that $A$ is algebraic.
Thus, the conclusion of Theorem~1 does not hold.
\end{itemize}

\section{Main Theorem in positive characteristic}

{\bf Lemma 9}.
Let $A$ be a finite-dimensional unital algebra over a field~$F$ of positive characteristic~$p$.
Then given an RB-operator $R$ of weight~$\lambda$ on $A$,
there exist $k,l\geq0$ such that $R^k(R+\lambda\id)^l = 0$.

{\sc Proof}.
\underline{Let $\lambda = 0$}. Since $A$ is finite-dimensional, there exists a minimal $l$ such that
\begin{equation}\label{ThmMain:p>0,lambda=0-Assum}
R^l(1) + \alpha_{l-1}R^{l-1}(1) + \ldots + \alpha_1 R(1) = 0,\quad \alpha_i\in F.
\end{equation}
It is well-known that unit does not lie in the image of $R$,
thus we take a zero free coefficient in~\eqref{ThmMain:p>0,lambda=0-Assum}.
Multiplying~\eqref{ThmMain:p>0,lambda=0-Assum} by $R(1)$
and applying the formula~\eqref{R^k(1)R(1)} for $\lambda=0$, we get
\begin{equation}\label{ThmMain:p>0,lambda=0-1}
(l+1)R^{l+1}(1) + l\alpha_{l-1}R^l(1) + \ldots + 2\alpha_1 R^2(1) = 0.
\end{equation}
Acting on~\eqref{ThmMain:p>0,lambda=0-Assum} by $R$, we have
\begin{equation}\label{ThmMain:p>0,lambda=0-2}
R^{l+1}(1) + \alpha_{l-1}R^l(1) + \ldots + \alpha_1 R^2(1) = 0.
\end{equation}
From~\eqref{ThmMain:p>0,lambda=0-1} and~\eqref{ThmMain:p>0,lambda=0-2}, we conclude
$$
\alpha_{l-1}R^l(1) + 2\alpha_{l-2}R^{l-1}(1) + \ldots + (l-1)\alpha_1 R^2(1) = 0.
$$
Since $l$ was chosen as a minimal number satisfying~\eqref{ThmMain:p>0,lambda=0-Assum},
we have two cases.

{\sc Case 1}.
The vectors $(1,\alpha_{l-1},\ldots,\alpha_2,\alpha_1)$
and $(\alpha_{l-1},2\alpha_{l-2}\ldots,(l-1)\alpha_1,0)$ are proportional
(with some nonzero coefficient).
It means that $\alpha_1 = \alpha_2 = \ldots = \alpha_{l-1} = 0$ and $R^l(1) = 0$, as required.

{\sc Case 2}.
The vector $(\alpha_{l-1},2\alpha_{l-2}\ldots,(l-1)\alpha_1,0)$ is zero.
If $l\leq p$, then we derive $\alpha_1 = \alpha_2 = \ldots = \alpha_{l-1} = 0$ and $R^l(1) = 0$.
If $l>p$, then only coefficients $\alpha_{l-p},\alpha_{l-2p},\ldots$ may be nonzero.
Anyway, we can find an equality
\begin{equation}\label{ThmMain:p>0,lambda=0-Fine}
R^{k_1}(1) + \beta_2 R^{k_2}(1) + \ldots + \beta_s R^{k_s}(1) = 0,\quad k_1>k_2>\ldots>k_s,\
p|k_1,\ldots,p|k_s
\end{equation}
with a minimal number of nonzero coefficients (here $s$).
If $s=1$, we are done.

Suppose that $s>1$. Consider $\Delta = k_1-k_s = p^tq$, where
$(p,q) = 1$ and $t\geq1$.
Acting by the operator $R^{\delta}$ on~\eqref{ThmMain:p>0,lambda=0-Fine}
for a number $\delta\in p\mathbb{N}$, we may get
$$
k_1 = p^N - p^t,\quad
k_s = p^N - p^t - \Delta
$$
for some sufficiently big~$N$.

Now, we mutiply~\eqref{ThmMain:p>0,lambda=0-Fine} by $R^{p^t}(1)$.
Applying Lemma~6b for $\lambda = 0$
$$
R^a(1)R^b(1) = \binom{a+b}{a}R^{a+b}(1),
$$
we have
$$
R^{k_1}(1)R^{p^t}(1) = 0,\quad
R^{k_s}(1)R^{p^t}(1) = \binom{p^t(p^{N-t}-q)}{p^t}R^{p^N - \Delta}(1),
$$
where the coefficient $\binom{p^t(p^{N-t}-q)}{p^t}$ is nonzero by the Lucas' Theorem.
So, we have found an expression like~\eqref{ThmMain:p>0,lambda=0-Fine}
with less number of nonzero coefficients. It is a contradiction, so $s = 1$, and we are done.

\underline{Let $\lambda$ be nonzero}, so we fix $\lambda = -1$.
Since $A$ is finite-dimensional, there exists an expression~$m$ of the form
$$
R^n(1) + \alpha_{n-1}R^{n-1}(1) + \ldots + \alpha_1 R^1(1) + \alpha_0 = 0, \quad \alpha_i\in F,
$$
of the minimal degree $n$.

Analogously to the proof of Lemma~7, we derive the formula
$$
q = (\alpha_{n-1}+n)R^n(1) + (2\alpha_{n-2}+(n-1)\alpha_{n-1})R^{n-1}(1) + \ldots + (n\alpha_0 + \alpha_1)R^1(1) = 0.
$$
If all $\alpha_i$ are zero, $R^n(1) = 0$ and we are done by Lemma~8.

Otherwise, find a minimal~$s$ such that $\alpha_s\neq0$.
So, $q = sm$, and we have the relations
\begin{equation}\label{ThmMain:p>0,system}
\begin{gathered} \allowdisplaybreaks
s\alpha_{s+1} = (n-s)\alpha_s + (s+1)\alpha_{s+1},\\
s\alpha_{s+2} = (n-s-1)\alpha_{s+1} + (s+2)\alpha_{s+2},\\
\ldots, \\
s\alpha_{n-1} = 2\alpha_{n-2} + (n-1)\alpha_{n-1},\\
s = \alpha_{n-1} + n.
\end{gathered}
\end{equation}
Define $n-s = tp + k$ for $0\leq k<p$ and $t\geq0$.
Due to the system~\eqref{ThmMain:p>0,system}, we have
\begin{multline*}
m = \alpha_s\left( R^s(1) - \binom{k}{1}R^{s+1}(1) + \binom{k}{2}R^{s+2}(1)
 \ldots + (-)^k \binom{k}{k}R^{s+k}(1) \right) \\
 + \alpha_{s+p}\left( R^{s+p}(1) -\ldots + (-)^k \binom{k}{k}R^{s+p+k}(1) \right) \\
 + \ldots
 + (-1)^k\left( R^{n-k}(1) - \ldots + (-1)^k \binom{k}{k}R^n(1) \right) \\
 = R^s(R-\id)^k ( \alpha_s + \alpha_{s+p} R^p + \ldots + (-1)^k R^{tp} )(1) = 0.
\end{multline*}
Acting on the last expression by the operator
$R^\delta(R-\id)^{\varepsilon}$
with suitable $\delta,\varepsilon$, we get the equality of the form
\begin{equation}\label{ThmMain:p>0,equ}
R^{pk_0} (R-\id)^p ( R^{pk_1} + \beta_2 R^{pk_2} + \ldots + \beta_e R^{pk_e})(1) = 0,
\end{equation}
where $k_1>k_2>\ldots>k_e$, $k_0>0$, and $\beta_i\neq0$, $i=2,\ldots,e$.

Since $(R-\id)^p = R^p -\id$ modulo $p$,
the equality~\eqref{ThmMain:p>0,equ} can be rewritten in the form
\begin{equation}\label{ThmMain:p>0,equ2}
R^{pm_1}(1) + \gamma_2 R^{p(m_1-m_2)}(1) + \ldots + \gamma_f R^{p(m_1-m_f)}(1) = 0,
\end{equation}
where $0<m_2<m_3<\ldots <m_f<m_1$.

We may assume, that $m_1 = p^{c-1}$ for some $c\geq1$.
Let $d$ be any natural not less than~$c$.
Acting by the operator $R^{p^d - p^c}$ on~\eqref{ThmMain:p>0,equ2},
we get the equation
\begin{equation}\label{ThmMain:p>0,equ3}
R^{p^d}(1) + \gamma_2 R^{p^d - pm_2}(1) + \ldots + \gamma_f R^{p^d -pm_f}(1) = 0.
\end{equation}

By Lemma~6b, we have for $\lambda = -1$
\begin{equation}\label{R^k(1)R^l(1)Nonzero}
R^{ap}(1)R^{bp}(1)
 =  \sum\limits_{c=0}^{\min\{a,b\}}(-1)^{cp}\binom{ap}{cp}\binom{p(a+b-c)}{ap}R^{p(a+b-c)}(1)
\end{equation}
over a field of characteristic~$p$.

With the help of the formula~\eqref{R^k(1)R^l(1)Nonzero} and the Lucas' Theorem,
we calculate
\begin{multline}
0 = (R^{p^d}(1) + \gamma_2 R^{p^d - pm_2}(1) + \ldots + \gamma_f R^{p^d -pm_f}(1))R^{p^d}(1) \\
 = E(d) := 2R^{2p^d}(1) + (-1)^{p^d}R^{p^d}(1) + \sum\limits_{j=2}^f\gamma_j \binom{2p^d-pm_j}{p^d} R^{2p^d - pm_j}(1).
\end{multline}
By the pigeonhole principle, we may find such $d_1>d_2\geq c$ that
$\binom{2p^{d_1}-pm_j}{p^{d_1}} = \binom{2p^{d_2}-pm_j}{p^{d_2}}$ for all $j=2,\ldots,f$.
Thus,
$$
R^{2p^{d_1}-2p^{d_2} }(E(d_2)) - E(d_1)
 = (-1)^p R^{p^{d_1} }( R^{ p^{d_1} - p^{d_2} } - \id )(1) = 0.
$$
Fix $\Delta = p^{d_1-1} - p^{d_2-1}$.
Then we have
\begin{equation}\label{ThmMain:p>0,equ4}
R^{p^k}(1) - R^{p^k-\Delta p}(1) = 0
\end{equation}
for all sufficiently big $k$.
Multiplying~\eqref{ThmMain:p>0,equ4} by $R^{p^k}(1)$,
we obtain
\begin{equation}\label{ThmMain:p>0,equ5}
2R^{2p^k}(1) + (-1)^p R^{p^k}(1) - \binom{2p^k-\Delta p}{p^k}R^{2p^k-\Delta p}(1) = 0.
\end{equation}
Acting by $\binom{2p^k-\Delta p}{p^k}R^{p^k}$ on~\eqref{ThmMain:p>0,equ4}
and subtracting the result from~\eqref{ThmMain:p>0,equ5}, we get
$$
\left(2 - \binom{2p^k-\Delta p}{p^k}\right)R^{2p^k}(1) + (-1)^p R^{p^k}(1) = 0.
$$
If the number in the brackets is zero, we are done.
Otherwise, we find an expression of the form
\begin{equation}\label{ThmMain:p>0,equ6}
R^{2p^k}(1) - lR^{p^k}(1) = 0,\quad l\in\mathbb{Z}^*_p.
\end{equation}
If $\charr(F) = p = 2$, then
$$
R^{2p^k}(1) - R^{p^k}(1)
 = R^{p^k}(R-\id)^{p^k}(1) = 0,
$$
and have proved the statement.

Suppose that $p>2$.
It is easy to show that consecutive applications of~\eqref{ThmMain:p>0,equ6}
imply the equality
\begin{equation}\label{ThmMain:p>0,equ7}
0
 = R^{p^{k+1} }(1) - l^{p-1}R^{p^k}(1)
 = R^{p^{k+1} }(1) - R^{p^k}(1) .
\end{equation}
Finally, with the help of~\eqref{ThmMain:p>0,equ7} we compute
\begin{multline*}
0 = (R^{p^{k+1} }(1) - R^{p^k}(1))^2 \\
  = 2R^{2p^{k+1} }(1) - R^{p^{k+1} }(1)
 - 2R^{p^{k+1} + p^k }(1)
 + 2R^{2p^k}(1) - R^{p^k}(1) \\
 = 2( R^{2p^k}(1) - R^{p^k}(1))
 = 2R^{p^k}(R-\id)^{p^k}(1),
\end{multline*}
so, $R^{p^k}(R-\id)^{p^k}(1) = 0$, and we are done. \hfill $\square$

{\bf Theorem 2}.
Let $A$ be a finite-dimensional unital algebra over
a field $F$ of positive characteristic. Let us fix $\lambda\in F$.

a) Given an RB-operator $R$ of weight~$\lambda$ on $A$,
there exist $k,l\geq0$ such that $R^k(R+\lambda\id)^l = 0$. In particular,
$\Spec(R)\subset\{0,-\lambda\}$.

b) There exists a natural~$n$ such that for every RB-operator~$R$
of weight~$\lambda$ on $A$ one can find $0\leq k\leq n$ such that
$R^k(R+\lambda\id)^{n-k}(1) = 0$ holds.
Moreover, $\rb_\lambda(A)\leq 2n$.

{\sc Proof}.
a) It is Lemma~9.

b) We may find such~$n$ by Lemma~9 and by the condition that $A$ is finite-dimensional.
For the bound $\rb_\lambda(A)\leq 2n$, we use Lemma~8 and repeat the end of the proof of Theorem~1.
\hfill $\square$

\section{Miminal polynomial of the image of unit}

Define polynomials
\begin{equation}\label{H-Poly}
H_{r,s}(x) = (x-r)(x-r+1)\ldots x(x+1)\ldots (x+s-1)
\end{equation}
for natural $r,s\geq0$.

{\bf Lemma 10}.
Let $A$ be a unital algebra and let $R$ be an RB-operator of weight~$-1$ on~$A$.
Denote $a = R(1)$. Then for all $r,s\geq0$,
\begin{equation}\label{ROnHrs}
(r+s+1)R(H_{r,s}(a)) = H_{r,s+1}(a).
\end{equation}

{\sc Proof}.
Let us prove~\eqref{ROnHrs} by induction on $r$.
For $r= 0$, we want to state the formula
\begin{equation}\label{R^n(a)}
(n+1)!R^n(a) = H_{0,n+1}(a)
\end{equation}
by induction on $n$.
For $n = 0$, we have $a = a$, it is true.
Suppose that we have proved that
$n!R^{n-1}(a) = H_{0,n}(a)$ or, equivalently,
$a^n = n!R^{n-1}(a) + f(a)$,
where
$f(x) = x^n - H_{0,n}(x)$.

By~\eqref{R^k(1)R(1)}, we calculate
\begin{multline*}
a^{n+1}
 = a^n\cdot a
 = (n!R^n(1)+f(a))R(1) \\
 = (n+1)!R^{n+1}(1) - n\cdot n!R^n(1) + af(a) \\
 = (n+1)!R^{n+1}(1) - nH_{0,n}(a) + a(a^n - H_{0,n}(a)) \\
 = (n+1)!R^{n+1}(1) + a^{n+1} - (a+n)H_{0,n}(a).
\end{multline*}
So, $(n+1)!R^{n+1}(1) = H_{0,n+1}(a)$ as required.
Thus, we get~\eqref{ROnHrs} for $r = 0$.

Suppose we have proved~\eqref{ROnHrs} for all numbers less than~$r$.
Represent
$$
H_{r,s}(x)
 = H_{r-1,s+1}(x) - (r+s)H_{r-1,s}(x).
$$
Applying the induction hypothesis for the last equality, we get
\begin{multline*}
(r+s+1)R(H_{r,s}(a))
 = (r+s+1)R(H_{r-1,s+1}(a)) - (r+s+1)(r+s)R(H_{r-1,s}(a)) \\
 = H_{r-1,s+2}(a) - (r+s+1)H_{r-1,s+1}(a)
 = H_{r,s+1}(a),
\end{multline*}
as required. \hfill $\square$

{\bf Lemma 11}.
Let $A$ be a unital algebra and let $R$ be an RB-operator of weight~$-1$ on~$A$.
Denote $a = R(1)$. Then for all $r,s\geq0$,
\begin{equation}\label{R'^rR^s}
(r+s)!(R-\id)^r R^s(1) = H_{r,s}(a).
\end{equation}

{\sc Proof}.
We prove the statement by induction on $r$.
For $r=0$, it follows from~\eqref{R^n(a)}.
Suppose we have proved~\eqref{R'^rR^s} for all numbers less than~$r$.
Applying the induction hypothesis and~\eqref{ROnHrs},
we get the following equalities
\begin{multline*}
(r+s)!(R-\id)^r R^s(1)
 = (R-\id)(r+s)H_{r-1,s}(a) \\
 = H_{r-1,s+1}(a) - (r+s)H_{r-1,s}(a)
 = H_{r-1,s}(a)(a+s-(r+s))
 = H_{r,s}(a).
\end{multline*}
So, Lemma is proved. \hfill $\square$

\newpage

{\bf Theorem 3}.
Let $A$ be a unital algebra over a field $F$
and let $R$ be an RB-operator of weight~$-1$ on $A$.
Denote $a = R(1)$.

a) If $\charr(F) = 0$ and $A$ is algebraic, then there exist $r,s$ such that $H_{r,s}(a) = 0$.

b) If $\charr(F) = p>0$, then $H_{0,p}(a) = 0$.

\noindent
In particular, we get that the minimal polynomial of $a$ under the conditions
has no multiple roots.

{\sc Proof}.
a) It follows from Lemma~7 and Lemma~11.

b) We are done by~\eqref{R'^rR^s}. \hfill $\square$

Now, we get easily the following result which requires in~\cite{Unital} several pages.

{\bf Corollary 1}~\cite[Thm.~4.17]{Unital}.
Given a field $F$ of characteristic~0, and an RB-operator of weight~$-1$ on $M_n(F)$,
the matrix $R(1)$ is diagonalizable and the set of diagonal elements in an appropriate basis
has the form $\{-s,-s+1,\ldots,0,1,\ldots,r-1,r\}$
for some natural numbers $r,s$ such that $r+s+1\leq n$.

We call an algebra $A$ as a nil-algebra if for every $x\in A$
there exists $k$ such that $x^k = 0$ for each bracketing of the nonassociative monomial $x^k$.
We generalize and simplify the proof of another statement from~\cite{Unital}.

{\bf Corollary 2}~\cite[Thm.~4.8]{Unital}.
Let $A$ be a unital algebra over a field~$F$ and
$A = F\cdot1\oplus N$ (direct vector-space sum), where $N$ is a nil-algebra.
Then each Rota---Baxter operator $R$ of nonzero weight on $A$ is splitting
and (up to $\phi$) we have $R(1) = 0$.

{\sc Proof}.
Denote $a = R(1)$. If $a\in N$, then $a^k = 0$ for some $k$.
By Theorem~3, we get $a = 0$. Then by Lemma~4, $R$ is splitting.

If $a\not\in N$, then $a = \alpha\cdot1 + n$, where $\alpha\in F$ and $n\in N$.
So, there exists $k$ such that $(a-\alpha\cdot1)^k = 0$.
By Theorem~3, $a-\alpha\cdot1 = 0$ or $a\in F$.
By Lemma~4, $R$ is splitting and $\phi(R)(1) = 0$.
\hfill $\square$

By~Corollary~2, we have $\rb_\lambda(\mathrm{Gr}_n) = 2$ for the Grassmann algebra~$\mathrm{Gr}_n$~\cite{Unital}.

{\bf Corollary 3}.
Let $A$ be a unital algebra over a field $F$ of characteristic zero
and let $R$ be an RB-operator of weight~$\lambda$ on~$A$.
Denote by $S$ a subalgebra generated by the set $\{R^k(1) \mid k\in\mathbb{N}\}$.

a) If $\dim(S) = \infty$, then $S$ is a free commutative RB-algebra of weight~$\lambda$ generated by~1.

b) If $\dim(S) = n <\infty$ and $\lambda = 0$, then
$S$~is isomorphic to the quotient of $F[x]$ by the ideal generated by $x^n$.
Moreover, $R$ acts on $S$ as follows, $R(x^i) = \frac{x^{i+1}}{i+1}$ for $i=0,\ldots,n-1$.

c) If $\dim(S)<\infty$, $\lambda = -1$, $F$ is algebraically closed,
and $H_{r,s}(x)$ is a minimal polynomial of $R(1)$, then $S\cong F^{r+s}$.

{\sc Proof}.
a) It follows directly.

b) It follows from Theorem~1 and formula~\eqref{Stirling1}.

c) Note that $S$ is isomoprhic to the quotient of $F[x]$ by the ideal generated
by the polynomial $H_{r,s}$ of degree $n = r + s$ which has no multiple roots.
So, $S$ is semisimple. Since $F$ is algebraically closed, $S\cong F^{r+s}$.
\hfill $\square$

The notion of capacity plays a crucial role in the structural theory of Jordan algebras~\cite{Taste}.
Let us generalize it for all unital algebras (of every variety).
Given a unital algebra~$A$, define a capacity of $A$ as the maximal number of pairwise
orthogonal idempotents in $A$. If capacity of $A$ equals~$n$, then there exist
$e_1,\ldots,e_n\in A$ such that $e_i^2 = e_i$ and $e_i e_j = 0$ if $i\neq j$.

{\bf Corollary 4}.
Let $A$ be a unital algebraic algebra over an algebraically closed field $F$ of characteristic zero.
If capacity of $A$ equals $n<\infty$, then $\rb_\lambda(A)\leq 2n$ for $\lambda\neq0$.

{\sc Proof}.
Let $R$ be an RB-operator of nonzero weight~$\lambda$ on $A$.
By Theorem~1a, there exist $k,l$ such that $R^k(R+\lambda\id)^l = 0$ on~$A$.
In particular, it means that the subalgebra~$S$ in~$A$ generated by the set
$\{R^j(1) \mid j\in\mathbb{N}\}$ is finite-dimensional.
By Corollary~3c, $S\cong F^m$ for some $m$.
Clearly, $m\leq n$ by the definition of capacity.
So, minimal polynomial of $P(1)$ for every RB-operator~$P$ of weight $\lambda$ on~$A$
has degree not greater than~$n$.
By the proof of Theorem~1b, $\rb_\lambda(A)\leq 2n$.
\hfill $\square$

\section{Rota---Baxter-index of different algebras}

\subsection{Sum of fields}
Let $F^n = Fe_1\oplus Fe_2\oplus \ldots \oplus Fe_n$
denote a direct sum of $n$ copies of a field $F$.
Here $e_i^2 = e_i$.

{\bf Proposition 1}.
$\rb_\lambda(F^n) = \begin{cases}
1, & \lambda = 0, \\
n, & \lambda \neq 0.
\end{cases}$

{\sc Proof}.
For $\lambda = 0$, it follows from Corollary~5.7~\cite{Unital}.

Let $\lambda\neq0$ and let $R$ be an RB-operator of weight~$\lambda$ on $F^n$.
By the dimensional reasons and the Main Theorem,
$R^k(R+\lambda\id)^{n-k} = 0$ for some $0\leq k\leq n$.
On the other hand, define $R$ on $F^n$ as follows:
$R(e_i) = e_{i+1}+\ldots+e_n$.
It is the RB-operator of weight~1 on~$F^n$.
Moreover, $R^n = 0$ and $R^{n-1}\neq0$. \hfill $\square$

\subsection{Matrix algebra}

Let us generalize the RB-operator (M1) on $M_2(F)$ from~\cite[Thm. 4.13]{Unital}
and the RB-operator 1-I on $M_3(F)$ from~\cite[Thm. 3]{GonGub}.

{\bf Example 1}.
A linear map $R$ defined on $M_n(F)$ as follows,
$$
R(e_{ij})
 = \begin{cases}
  \sum\limits_{k\geq1}e_{i+k,j+k}, & i\geq j, \\
- \sum\limits_{k\geq0}e_{i-k,j-k}, & i<j,
   \end{cases}
$$
is an RB-operator of weight~1 on $M_n(F)$.
Moreover, $R^n(R+\id)^{n-1} = 0$ but we have
$(R(R+\id))^{n-1}\neq0$ and $R^n(R+\id)^{n-2}\neq0$.

The following theorem was stated in~\cite{Unital} in the case of $\lambda = 0$.

{\bf Theorem 4}.
Let $F$ be a field of characteristic zero. Then
$\rb_\lambda(M_n(F)) = 2n-1$ for all $\lambda\in F$.

{\sc Proof}.
For $\lambda = 0$, it was proved in~\cite{Unital}.

Let $\lambda\neq0$, we may assume that $\lambda = 1$.
By Example~1, $\rb_\lambda(M_n(F)) \geq 2n-1$.

Let $R$ be an RB-operator of weight~$-1$ on $M_n(F)$.
By the Main Theorem, $\Spec(R)\subset\{0,1\}$.
So, we extend the action of $R$ to the algebra $M_n(\bar{F})$.
We want to show that either $R^n(R-\id)^{n-1} = 0$ or $R^{n-1}(R-\id)^n = 0$.
Thus, it is enough to state this for $R$ over $\bar{F}$.

Define $a = R(1)$.
If the degree of the minimal polynomial $m_t(a)$ of~$a$ is less than~$n$, then
by Theorem~1, $(R(R-\id))^{n-1} = 0$, and we are done.
Suppose that $\deg(m_t(a)) = n$ and $m_t(a)$ coincides with
the characteristic polynomial of $a$.
By Corollary~1, $a$ is conjugate to a diagonal matrix
$D = \diag\{-r,-r+1,\ldots,s-2,s-1\}$ with $r + s = n$.
By Lemma~1 from \cite{GonGub}, all subspaces
$V_q = \Span\{e_{ij}\mid i-j = q\}$ are $R$-invariant.
Since all $V_q$ for $q\neq0$ have dimension less or equal to $n-1$,
we have $(R(R-\id))^{n-1} = 0$ on each of such $V_q$.
Further, $V_0$ which is the subalgebra of diagonal matrices from $M_n(\bar{F})$
has the dimension equal to $n$. So, either $R^n = 0$ or $(R-\id)^n$ or
$(R(R-\id))^{n-1} = 0$ on~$V_0$. In any case, $\rb_\lambda(M_n(\bar{F})) = 2n-1$
as well as of $M_n(F)$. \hfill $\square$

Recall that $t_n(F)$ denotes the subalgebra of not strictly upper-triangular matrices from $M_n(F)$.

{\bf Corollary 5}.
Let $F$ be a field of characteristic zero and let $\lambda\in F$ be nonzero.
Then $\rb_\lambda(t_n(F)) = 2n-1$.

{\sc Proof}.
For simplicity, let $\lambda = 1$.
The restriction of the RB-operator $R$ from Example~1 on $t_n(F)$ implies
$\rb_\lambda(t_n(F)) \geq 2n-1$.
We may assume that $F$ is algebraically closed, since
$\rb_\lambda(t_n(F))\leq \rb_\lambda(t_n(\bar{F}))$.
If the degree of the minimal polynomial $m_x(a)$ of $a = R(1)$
is less than $n$, we are done by Theorem~1.
Suppose that $\deg(m_x(a)) = n$.

Find $\psi\in\Aut(M_n(F))$ such that
$R^{(\psi)}(1)$ equals to the $D = \diag\{-r,-r+1,\ldots,s-2,s-1\}$ with $r + s = n$.
Now, the RB-operator $P = R^{(\psi)}$ acts on the subalgebra $\psi(t_n(F))$.
Define $V_q = \Span\{e_{ij}\mid i-j = q\}$.
Applying the proof of Lemma~1 from \cite{GonGub},
we conclude that a subspace $V_q\cap \psi(t_n(F))$ is $P$-invariant for every~$q$.
It remains to repeat the proof of Theorem~4. \hfill $\square$

{\bf Conjecture 1}.
Let $F$ be a field of characteristic zero. Then $\rb_0(t_n(F)) = n$.

Let us extend Conjecture 5.24 from~\cite{Unital} on the case of nonzero weight.

{\bf Conjecture 2}.
Let $A = \bigoplus\limits_{i=1}^k M_{n_i}(F)$
be a semisimple finite-dimensional associative algebra
over a field $F$ of characteristic zero. Then

a) $\rb_0(A) = 2\max\limits_{i=1,\ldots,k}\{n_i\} - 1$.

b) $\rb_\lambda(A) = \sum\limits_{i=1}^k n_i + \max\limits_{i=1,\ldots,k}\{n_i\} - 1$ for nonzero~$\lambda$.

The following example shows that there exists an RB-operator reaching the bound from Conjecture~2b.

Let $D_n$ denote the subalgebra of all diagonal matrices in $M_n(F)$
and $L_n$ ($U_n$) the set of all strictly lower (upper) triangular
matrices in $M_n(F)$.

{\bf Example 2}.
We define a triangular-splitting RB-operator $R$ of nonzero weight~$\lambda$ on
the algebra $A = \bigoplus\limits_{i=1}^k M_{n_i}(F)$ with the subalgebras
$A_0 = D$, $A_- = U$, $A_+ = L$, where
$$
D = \bigoplus\limits_{i=1}^k D_{n_i},\quad
L = \bigoplus\limits_{i=1}^k L_{n_i},\quad
U = \bigoplus\limits_{i=1}^k U_{n_i},
$$
and $R$ is defined on $D \cong F^N$ for $N = \sum\limits_{i=1}^k n_i$
in the same way as the RB-operator from Example~1.
So, $x^N(x+\lambda)^{n_0-1} = 0$ is the minimal polynomial for~$R$.
Here $n_0 = \max\limits_{i=1,\ldots,k}\{n_i\}$.

\subsection{Simple Jordan algebra of a bilinear form}
Let $J = J_{n+1}(f) = F1\oplus V$ be a direct vector-space sum of $F$
and finite-dimensional vector space~$V$, $\dim V = n > 1$,
and $f$ be a nondegenerate symmetric bilinear form on~$V$.
Under the product
\begin{equation}\label{FormProduct}
(\alpha\cdot1 + a)(\beta\cdot1 + b)
 = (\alpha\beta+f(a,b))\cdot1 + (\alpha b +\beta a),\quad
\alpha,\beta\in F,\ a,b\in V,
\end{equation}
the space $J$ is a simple Jordan algebra \cite{Nearly}.
We assume that $\charr F\neq2,3$.

Let us choose a basis $e_1$, $e_2$, \ldots, $e_n$ of $V$ such that
the matrix of the form $f$ in this basis is diagonal with
elements $d_1,d_2,\ldots,d_n$ on the main diagonal.
As $f$ is nondegenerate, $d_i\neq0$ for each $i$.

{\bf Example 3}~\cite{BGP}.
Let $J_{2n}(f)$, $n\geq2$, be the simple Jordan algebra of bilinear form~$f$
over an algebraically closed field~$F$ with $\charr(F)\neq2,3$.
Let $R$ be a linear operator on $J_{2n}(f)$
defined by a matrix $(r_{ij})_{i,j=0}^{2n-1}$
in the basis $1,e_1,e_2,\ldots,e_n$
with the following nonzero entries
\begin{gather*}
r_{00} = -3,\quad r_{01} = \sqrt{d_1},\quad r_{10}
 = -\frac{1}{\sqrt{d_1}},\quad r_{jj} = -1,\ j=1,\ldots,2n-1, \\
r_{i\,i+1} = \frac{d_{i+1}}{d_i}\sqrt{ -\frac{d_i}{d_{i+1}} }, \quad
r_{i+1\,i} = - \sqrt{ -\frac{d_i}{d_{i+1}} },\ i=2,\ldots,2n-2.
\end{gather*}
Then $R$ is an RB-operator of weight~$2$ on $J_{2n}(f)$.
Moreover, $R(R+2\id)^2 = 0$ but $R(R+2\id)\neq0$ and $(R+2\id)^2 \neq 0$.

{\bf Proposition 2}.
Let $\lambda\in F$ be nonzero and let $F$ be an algebraically closed field
with $\charr F\neq2,3$. Then
$\rb_\lambda(J_n(F))
 = \begin{cases}
2, & n \mbox{ is odd}, \\
3, & n \mbox{ is even}.
\end{cases}$

{\sc Proof}.
When $n$ is odd, the statement follows from \cite[Thm.\,4.1]{BGP}.

Suppose that $n$ is even. By Example~3, we have $\rb_\lambda(J_n(F))\geq3$.
Since $F$ is algebraically closed, we may assume that all $d_i$ equal to 1.
Let $R$ be an RB-operator of weight~$\lambda$ on $J_n(F)$.
Assume that $R$ is defined by a matrix $(r_{ij})_{i,j=0}^n$
in the basis $1,e_1,e_2,\ldots,e_n$.
Due to the proof of Theorem~4.1 from~\cite{BGP}, we have
the following relations
\begin{gather*}
r_{ii} = -\frac{\lambda}{2},\quad
r_{kl} = - r_{lk},\quad
r_{0i} = zr_{i0},\\
\sum\limits_{p=1}^n \bar{r}_{p0}^2 = \frac{\lambda^2}{4},\quad
\sum\limits_{p=1}^n \bar{r}_{pi}^2 = \bar{r}_{0i}^2,\quad
\sum\limits_{p=1}^n \bar{r}_{pk}\bar{r}_{pl}
 = \bar{r}_{0k}\bar{r}_{0l},\quad
\sum\limits_{p=1}^n \bar{r}_{pi}\bar{r}_{p0}
 = -\frac{\lambda z}{2}\bar{r}_{0i},
\end{gather*}
where $i,k,l>0$ and $k\neq l$,
and we have three subcases:

(I) $z = 1$, $\bar{r}_{00} = -\lambda/2$,

(II) $z = -1$, $\bar{r}_{00} = \lambda/2$,

(III) $z = -1$, $\bar{r}_{00} = -3\lambda/2$.

It is easy to check that an RB-operator $R$ from (I) is splitting.
Further, $\phi$ maps RB-operators from (II) to (III) and vice versa.
Finally, an RB-operator from (III) satisfies the equation
$R(R+\lambda\id)^2 = 0$. \hfill $\square$

\subsection{Matrix Cayley---Dickson algebra}

Cayley---Dickson algebras are known to be only simple nonassociative alternative
finite-dimensional algebras, they are either division algebras or so called split ones.
The first ones have only trivial RB-operators~\cite{BGP}.
Let us study the $\rb(\lambda)$-index of the matrix (split)
Cayley---Dickson algebra $C(F)$, which is simple alternative 8-dimensional algebra.

Recall the definition of the product on $C(F) = M_2(F)\oplus v M_2(F)$.
We extend the product from
$M_2(F)$ on $C(F)$ as follows:
\begin{equation}\label{CD-product}
a(vb) = v(\bar{a}b),\quad (vb)a = v(ab),\quad
(va)(vb) = b\bar{a},\quad a,b\in M_2(F).
\end{equation}
Here\, $\bar{}$\, denotes the involution on $M_2(F)$
which maps a matrix $a = (a_{ij})\in M_2(F)$ to the matrix
$\begin{pmatrix}
a_{22} & - a_{12} \\
-a_{21} & a_{11} \end{pmatrix}$.

Given an algebra $A$ with the product $\cdot$, define the operations $\circ$, $[,]$ on
the space $A$:
$$
a\circ b=a\cdot b+b\cdot a,\quad [a,b]=a\cdot b-b\cdot a.
$$
We denote the space $A$ endowed with $\circ$ as $A^{(+)}$
and the space $A$ endowed with $[,]$ as~$A^{(-)}$.

{\bf Lemma 12}.
Let $A$ be an algebra. Then

(a)~\cite{BGP}
If $R$ is an RB-operator of weight~$\lambda$ on~$A$,
then $R$ is an RB-operator of weight~$\lambda$ on $A^{(\pm)}$.

(b) $\rb_\lambda(A)\leq \rb_\lambda(A^{\pm})$.

{\sc Proof}.
(b) It follows from the definition of $\rb(\lambda)$-index of an algebra and (a).
\hfill $\square$

{\bf Proposition 3}.
Let $F$ be an algebraically closed field with $\charr(F)\neq2,3$.
Then $\rb_\lambda(C(F)) = 3$ for all $\lambda\in F$.

{\sc Proof}.
When $\lambda = 0$, it was proved in \cite[Thm.\,5.26]{Unital}.

Let $\lambda\neq0$. It is known~\cite[p.\,57]{Nearly} that
$C(F)^{(+)}$ is isomorphic to the simple Jordan algebra of
a~bilinear form. By Proposition~2 and Lemma~12b, $\rb_\lambda(C(F))\leq3$.

Consider a triangular-splitting RB-operator $R$ of weight~1 on $C(F)$
with subalgebras $A_0 = \Span\{e_{11},e_{22},ve_{11},ve_{22}\}$,
$A_+ = \Span\{e_{12},ve_{12}\}$, and $A_- = \Span\{e_{21},ve_{21}\}$.
Further, $R$ is defined on $A_0$ again as a triangular-splitting
RB-operator with subalgebras
$B_0 = \Span\{e_{11},e_{22}\}$,
$B_+ = \Span\{ve_{11}\}$, and
$B_- = \Span\{ve_{22}\}$. Finally, $R$ is defined on $B_0\cong F\oplus F$
as follows: $R(e_{11}) = e_{22}$, $R(e_{22}) = 0$.
Thus, $R^2(R+\id) = 0$ on the entire $C(F)$
but $R^2\neq0$ and $R(R+\id)\neq0$.
So, $\rb_\lambda(C(F)) = 3$. \hfill $\square$

\section{Applications}

\subsection{Rota---Baxter modules}

Let $A$ be an algebra and let $R$ be an RB-operator $R$ of weight~$\lambda$ on $A$.
An $A$-bimo\-dule~$M$ with a linear map $p\colon M\to M$ is called
a~Rota---Baxter bimodule over $(A,R)$ and simply an $(A,R)$-bimodule~\cite{RBModules},
if the following equalities hold for all $x\in A$ and for all $m\in M$,
\begin{equation}\label{RBModule}
\begin{gathered}
R(x)p(m) =  p( R(x)m + xp(m) + \lambda xm ),\\
p(m)R(x) =  p( mR(x) + p(m)x + \lambda mx ).
\end{gathered}
\end{equation}
We call an $(A,R)$-bimodule $(M,p)$ a unital one, if
$M$ is a unital $A$-bimodule.

The relations~\eqref{RBModule} are equivalent to the condition that
$R+p$ is an RB-operator of weight~$\lambda$ on $A\ltimes M$.

{\bf Proposition 4}.
Let $A$ be a unital algebraic algebra over a field $F$ of characteristic zero
and let $R$ be an RB-operator of weight~$\lambda$ on $A$.
Also, $(M,p)$ is a unital $(A,R)$-bimo\-dule.
Then

a) If $\lambda = 0$, then $p$ is nilpotent on $M$.

b) If $\lambda\neq0$, then there exist $s,t$ such that
$p^s(p+\lambda\id)^t = 0$ on $M$. So, $M = M_0\oplus M_1$,
where $M_0 = \ker(p^s)$ and $M_1 = \ker(p+\lambda\id)^t$.
For $k,l$ such that $R^k(R+\lambda\id)^l = 0$ define
$A_0 = \ker(R^k)$ and $A_1 = \ker(R+\lambda\id)^l$.
We get that $(M_i,p)$ is a (not necessary unital) $(A_i,R)$-bimodule for $i=0,1$.

c) $\rb_\lambda(A\ltimes M)\leq 2\rb_\lambda(A)$ for all $\lambda\in F$.

{\sc Proof}.
a) By Theorem~1a, $R+p$ is nilpotent on $A\ltimes M$, so $p$ is nilpotent on $M$.

b) By Theorem~1a, there exist $n,m$ such that $(R+p)^n (R+p+\lambda\id)^m = 0$ on $A$.
So, $p^n(p+\lambda\id)^m = 0$ on $M$, and
$M = M_0\oplus M_1$ for $M_0 = \ker(p^s)$ and $M_1 = \ker((p+\lambda\id)^t )$.
Also, $\ker(R+p)^n$ and $\ker(R+p+\lambda\id)^m$ are subalgebras of $A\ltimes M$.
It implies that $(M_0,p)$ is an $(A_0,R)$-bimodule and $(M_1,p)$ is an $(A_1,R)$-bimodule.

c) Suppose that $\rb_\lambda(A)<\infty$.
Then there exists a minimal $n$ such that for every RB-operator
$R$ of weight~$\lambda$ on $A$ one can find
$0\leq k\leq n$ satisfying the equality
$R^k(R+\lambda\id)^{n-k}(1) = 0$.
On the one hand, $n\leq\rb_\lambda(A)$ by the definition of the $\RB(\lambda)$-index.
On the other hand, since
$(R+p)^k(R+\lambda\id)^{n-k}(1) = R^k(R+\lambda\id)^{n-k}(1) = 0$,
we get $\rb_\lambda(A)\leq 2n$ by Theorem~1b.
\hfill $\square$

{\bf Remark 2}.
The same result holds if $A$ is a finite-dimensional unital
algebra over any field (not necessary of characterstic zero).

\subsection{Parameterized Rota---Baxter operators}

The study of the parameterized associative Yang---Baxter equation
was initiated by A.~Polishchuk~\cite{Polishchuk} in 2002.
Further, the definition of parameterized associative Yang---Baxter equation
of weight $\lambda$ appeared~\cite{FardThesis,Anghel}.
Given an associative algebra~$A$, a linear map
$r\colon F\to A\otimes A$ is called a solution of
the parameterized associative Yang---Baxter equation
of weight~$\lambda$ if the equation
\begin{equation}\label{AYBE}
r_{13}(z_{13})r_{12}(z_{12})
 - r_{12}(z_{12})r_{23}(z_{23})
 + r_{23}(z_{23})r_{13}(z_{13})
 + \lambda r_{13}(z_{13}) = 0
\end{equation}
holds for all $z_1,z_2,z_3\in F$.
Here $z_{ij} = z_i-z_j$ and given a decomposition
$r(z) = \sum\limits_i a_i\otimes b_i$,
$$
r_{12} = \sum a_i\otimes b_i\otimes 1,\quad
r_{13} = \sum a_i\otimes 1\otimes b_i,\quad
r_{23} = \sum 1\otimes a_i\otimes b_i.
$$

It is known~\cite{FardThesis} that a solution $r$ of~\eqref{AYBE}
gives rise to a parameterized Rota---Baxter operator $R$ of weight~$\lambda$ on $A$.
More precisely, a linear map
$R\colon F\otimes A\to A$ is called parameterized Rota---Baxter operator $R$
of weight~$\lambda$ on $A$ if $R$ satisfies the relation
\begin{equation}\label{RBparam}
R_b(x)R_c(y) = R_c( R_{b-c}(x)y ) + R_b( xR_{c-b}(y) )  + \lambda R_b(xy),
\end{equation}
where $b,c\in F$ and $x,y\in A$.
By $R_b(x)$ we denote the element $R(b,x)\in A$.

So, a parameterized RB-operator $R$ on $A$ is a family of linear operators
$R_b\colon A\to A$, $b\in F$.
Note that $R_0$ is a usual RB-operator of weight~$\lambda$ on $A$.

Let us apply the results of Theorem~1 to get some information about operators $R_b$.

{\bf Proposition 5}.
Let $A$ be a unital algebraic associative algebra over a field of characteristic zero.
Let $R$ be a parameterized RB-operator of weight~$\lambda$ on $A$. Then

a) If $\lambda = 0$, then there exists $n$ such that
$(R_b(1))^n = 0$ for all $b\in F$,

b) If $\lambda  = -1$, then there exist $n,m$
such that $H_{n,m}(R_b(1)) = 0$ for all $b\in F$.

{\sc Proof}.
a) By Theorem~1a, there exists $n$ such that $R^n_0(1) = 0$.
By~\eqref{RBparam}, we get $(R_b(1))^{n+1} = R_b(R_0^n(1)) = 0$.

b) Denote $a = R_0(1)$ and $c = R_b(1)$.
Note that by~\eqref{RBparam},
$$
R_b(x)R_b(y) = R_b(R_0(x)y + xR_0(y) - xy).
$$
So, the right-hand side has the form
$R_b(Z)$, where $Z$ is a linear combination of the products
involving only the action of $R_0$.
The same holds for the product of every number of elements $R_b(x_i)$.
This fact allows us analogously to the formula
$$
(R_0(1))^n
 = a^n
 = R_0\left( \sum\limits_{k=1}^n k! S(n,k)(-1)^{n-k} R_0^{k-1}(1) \right),
$$
from~\cite[\S3.3.2]{GuoMonograph} to state the fomula
$$
(R_b(1))^n
 = c^n
 = R_b\left( \sum\limits_{k=1}^n k! S(n,k)(-1)^{n-k} R_0^{k-1}(1) \right).
$$
Here $S(n,k)$ denotes a Stirling number of the second kind.

So, if we have an equility
$g(a) = R_0(q(a))$ for some polynomials $g,q\in F[x]$,
then we may derive the equality $g(c) = R_b(q(a))$.
In particular, applying Lemma~10, we get
$$
\frac{H_{r,s+1}(a) }{r+s+1}
 = R_0( H_{r,s}(a) ),\quad
\frac{H_{r,s+1}(c) }{r+s+1}
 = R_b( H_{r,s}(a) ).
$$
By Theorem~1a and Lemma~11, there exist $r,s$ such that $H_{r,s}(a) = 0$.
Thus, $H_{r,s+1}(c) = 0$, we are done.
\hfill $\square$

Suppose that $R$ is a parameterized RB-operator of weight~$\lambda$ on $M_n(\mathbb{C})$.
So, Proposition~5 says that if $\lambda = 0$, then
matrices $R_b(1)$ are nilpotent for all $b\in F$.
If $\lambda \neq0$, then matrices $R_b(1)$ are diagonalizable for all $b\in F$.

\subsection{Derived Rota---Baxter algebras}

Let us generalize the following construction already known
for different varieties of algebras, see~\cite{GuoMonograph,preLie,Burde}.

Let $A$ be an algebra of a variety $\Var$ and let $R$ be an RB-operator of weight~$\lambda$ on $A$.
Denote by $V$ the vector space $A$.
We define by induction the algebra $A_i$, $i\geq0$, on $V$ as follows
\begin{equation}\label{DerivedProduct}
\begin{gathered}
x\circ_0 y = xy,\\
x\circ_{i+1} y = R(x)\circ_i y + x\circ_i R(y) + \lambda x\circ_i y.
\end{gathered}
\end{equation}
It is known~\cite{BBGN2011,Embedding} that $R$ is an RB-operator
of weight~$\lambda$ on all $A_i$ and all algebras $A_i$ belong to the same variety $\Var$.
By the definition, $R$ and $-(R+\lambda\id)$ are homomorphisms from $A_n$ to $A_{n-1}$,
and we get compositions of such homomorphisms,
$$
\ldots \xrightarrow[-R-\lambda \id]{R} A_{i+1} \xrightarrow[-R-\lambda \id]{R} A_i \xrightarrow[-R-\lambda \id]{R} A_{i-1}
 \xrightarrow[-R-\lambda \id]{R} \cdots  \xrightarrow[-R-\lambda \id]{R} A_0.
$$

Let us call $A_i$ as the $i$th derived Rota---Baxter algebra of the RB-algebra $A$.
If there exists a~minimal natural number $m$ such that $A_n\cong A_m$ for all $n>m$,
then we define the limiting derived algebra $A_\infty$ as $A_m$.
Otherwise, $A_\infty$ is not defined.

The following result generalizes Corollary~4.10 from~\cite{preLie} in the case of weight zero.

{\bf Proposition 6}.
Let $A$ be a unital algebraic algebra over a field $F$ of characteristic zero.
Let $R$ be an RB-operator of weight~$\lambda$ on $A$. Then $A_\infty$ is well-defined.

a) If $\lambda = 0$, then $A_\infty$ has trivial product.

b) If $\lambda\neq0$, then
$A_\infty \cong \ker(R^n)\oplus \ker(R+\lambda\id)^n$ for sufficiently big $n$.

{\sc Proof}.
a) By Theorem~1, we know that $R^n = 0$ for some natural $n$.
By~\eqref{DerivedProduct}, we get
\begin{equation}\label{DerivedProductZero}
x\circ_k y = \sum\limits_{i=0}^k \binom{k}{i}R^i(x)R^{k-i}(y).
\end{equation}
Thus, $x\circ_i y = 0$ for all $i\geq2n-1$.

b) By Theorem~1, there exists $n$ such that $R^n(R+\lambda\id)^n = 0$.
So, $\ker(R^n)$ and $\ker(R+\lambda\id)^n$ are ideals in $A_t$ for all $t\geq n$
as kernels of the homomorphisms.
So, $A_t = \ker(R^n)\oplus \ker(R+\lambda\id)^n$ for all $t\geq n$.
It remains to note that $R$ is invertible on $\ker(R+\lambda\id)^n$
as well as $R+\lambda\id$ on $\ker(R^n)$.
Thus,
$$
\langle \ker(R^n),\circ_s\rangle \cong \langle \ker(R^n),\circ_n\rangle,\quad
\langle \ker(R+\lambda\id)^n,\circ_s\rangle \cong \langle \ker(R+\lambda\id)^n,\circ_n\rangle
$$
for all $s\geq0$. It implies that $A_t\cong A_n$ for all $t\geq n$.
\hfill $\square$

{\bf Remark 3}.
The same result holds if $A$ is a finite-dimensional unital
algebra over any field (not necessary of characterstic zero).

{\bf Question}.
Let $A = M_n(F)$. What is a minimal $m$ such that $A_m$
has trivial product for every RB-operator $R$ of weight zero on $A$?

\section*{Acknowledgements}

Author is grateful to the participants of the scientific seminar
in Altai State Pedagogical University (Barnaul)
for the helpful discussion.
Author is thankful to P. Kolesnikov for the helpful remarks.

Author is supported by the Program of fundamental scientific researches
of the Siberian Branch of Russian Academy of Sciences, I.1.1, project 0314-2019-0001.

\noindent Vsevolod Gubarev \\
Sobolev Institute of Mathematics \\
Acad. Koptyug ave. 4, 630090 Novosibirsk, Russia \\
Novosibirsk State University \\
Pirogova str. 2, 630090 Novosibirsk, Russia \\
e-mail: wsewolod89@gmail.com
\end{document}